\documentclass[12pt]{article}
\usepackage{amsthm,amsfonts,amssymb,amscd}

\usepackage{amsmath}
\usepackage{color}
\usepackage{enumerate}
\usepackage{graphicx}
\usepackage{epstopdf}

\begin{document}

\begin{center}
\large \bf Birationally rigid complete intersections\\
of codimension two
\end{center}\vspace{0.3cm}

\centerline{Daniel Evans and Aleksandr Pukhlikov}\vspace{0.3cm}

\parshape=1
3cm 10cm \noindent {\small \quad\quad\quad \quad\quad\quad\quad
\quad\quad\quad {\bf }\newline We prove that in the parameter
space of $M$-dimensional Fano complete intersections of index one
and codimension two the locus of varieties that are not
birationally superrigid has codimension at least $\frac12
(M-9)(M-10)-1$.

Bibliography: 25 titles.} \vspace{0.3cm}

\vspace{0.3cm}

\section*{Introduction}

{\bf 0.1. Statement of the main result.} Birational (super)
rigidity is known for almost all families of Fano complete
intersections of index one in the projective space, see
\cite{Pukh13z,Pukh14a,Pukh13a}. Typically birational superrigidity
was shown for a generic (in particular, non-singular) variety in
the family. Now the improved techniques make it possible to obtain
more precise results, covering complete intersections with certain
simple types of singularities and estimating the codimension of
the subset of non-rigid varieties in the parameter space of the
family. The first work of this type for a family of Fano varieties
was done in \cite{EP} for Fano hypersurfaces of index 1. Here we
do it for complete intersections of codimension two.\vspace{0.1cm}

In this paper, the symbol ${\mathbb P}$ stands for the complex
projective space ${\mathbb P}^{M+2}$, where $M\geq 13$. Fix two
integers $d_2\geq d_1\geq 2$, such that $d_1+d_2=M+2$ and consider
the space
$$
{\cal P}={\cal P}_{d_1,M+3}\times{\cal P}_{d_2,M+3}
$$
of pairs of homogeneous polynomials $(f_1,f_2)$ on ${\mathbb P}$
(that is to say, in $M+3$ variables
$x_0,\dots,x_{M+2}$\vspace{0.1cm}) of degrees $d_1$ and $d_2$,
respectively. The symbol $V(f_1,f_2)$ denotes the set of common
zeros of $f_1$ and $f_2$. The following claim is the main result
of this paper.\vspace{0.1cm}

{\bf Theorem 0.1.} {\it There exists a Zariski open subset ${\cal
P}_{\rm reg}\subset{\cal P}$ such that:\vspace{0.1cm}

{\rm (i)} for every pair $(f_1,f_2)\in{\cal P}_{\rm reg}$ the
closed set $V=V(f_1,f_2)$ is irreducible, reduced and of
codimension 2 in ${\mathbb P}$ with the singular locus
$\mathop{\rm Sing}V$ of codimension at least 7 in $V$, so that $V$
is a factorial projective algebraic variety; the singularities of
$V$ are terminal, so that $V$ is a primitive Fano variety of index
1 and dimension $M$;\vspace{0.1cm}

{\rm (ii)} the estimate
$$
\mathop{\rm codim} (({\cal P}\backslash{\cal P}_{\rm
reg})\subset{\cal P}) \geq\frac12(M-9)(M-10)-1
$$
holds;\vspace{0.1cm}

{\rm (iii)} for every pair $(f_1,f_2)\in{\cal P}_{\rm reg}$ the
Fano variety $V=V(f_1,f_2)$ is birationally
superrigid.}\vspace{0.1cm}

See \cite[Chapter 2]{Pukh13z} for the definitions of birational
rigidity and superrigidity as well as for the standard
implications of these properties: Theorem 0.1 implies that for
every pair $(f_1,f_2)\in{\cal P}_{\rm reg}$ the corresponding Fano
complete intersection $V=V(f_1,f_2)\subset{\mathbb P}$ admits no
structures of a rationally connected fibre space, that is to say,
there exists no rational dominant map $\varphi\colon
V\dashrightarrow S$ onto a positive dimensional base $S$ such that
the fibre of general position is rationally connected. In
particular, $V$ is non-rational. Another well known implication is
that the groups of birational and biregular self-maps of $V$ are
the same: $\mathop{\rm Bir}V=\mathop{\rm Aut} V$.\vspace{0.1cm}

Now we describe the set ${\cal P}_{\rm reg}$ by explicit
conditions (some of them are global but most of them are local)
and outline the proof of Theorem 0.1.\vspace{0.3cm}


{\bf 0.2. Regular complete intersections.} Consider a pair of
homogeneous polynomials $(f_1,f_2)\in{\cal P}$, both non-zero.
Below we list the conditions that these polynomials are supposed
to satisfy for a regular pair.\vspace{0.1cm}

(R0.1) The polynomial $f_1$ is irreducible and the hypersurface
$\{f_1=0\}=F_1$ has at most quadratic singularities of rank
5.\vspace{0.1cm}

{\bf Remark 0.1.} This condition ensures that $F_1$ is a factorial
variety so that $\mathop{\rm Cl} F_1\cong\mathop{\rm Pic} F_1$ is
generated by the class of a hyperplane section and every effective
divisor on $F_1$ is cut out by a hypersurface in ${\mathbb
P}$.\vspace{0.1cm}

(R0.2) $f_2|_{F_1}\not\equiv 0$ and moreover the closed set
$\{f_2|_{F_1}=0\}$ is irreducible and reduced.\vspace{0.1cm}

(R0.3) Every point $o\in V=V(f_1,f_2)$ is

\begin{itemize}

\item either non-singular,

\item or a quadratic singularity,

\item or a biquadratic singularity.

\end{itemize}

For each of the three types the local regularity conditions will
be stated separately. Given a point $o\in V$, we fix a system of
affine coordinates $z_1,\dots,z_{M+2}$ on an affine subset $o\in
{\mathbb A}^{M+2}\subset{\mathbb P}^{M+2}$ with the origin at $o$,
and write down the expansions of the polynomials $f_i$:
$$
\begin{array}{ccccccccccccc}
f_1 & = & q_{1,1} & + & q_{1,2} & + & \dots & + & q_{1,d_1}, & & & &\\
f_2 & = & q_{2,1}& +  & q_{2,2} & + & \dots & + & q_{2,d_1} & + &
\dots &
+ & q_{2,d_2},\\
\end{array}
$$
where $q_{i,j}$ are homogeneous of degree $j$. We list the
homogeneous polynomials in the {\it standard order} as follows:
$$
q_{1,1},q_{2,1},q_{1,2},q_{2,2},\dots,q_{1,d_1},q_{2,d_1},\dots,q_{2,d_2},
$$
so that polynomials of smaller degrees precede the polynomials of
higher degrees and for $j\leq d_1$ the form $q_{1,j}$ precedes
$q_{2,j}$.\vspace{0.1cm}

Every non-singular point $o\in V$ is assumed to satisfy the
regularity condition\vspace{0.1cm}

(R1) the polynomials $q_{i,j}$ in the standard order with the last
two of them removed form a regular sequence in ${\cal
O}_{o,{\mathbb P}}$.\vspace{0.1cm}

Every quadratic point $o\in V$ is assumed to satisfy a number of
regularity conditions. Note that in this case at least one of the
linear forms $q_{1,1},q_{2,1}$ is non-zero and the other one is
proportional to it. We denote a non-zero form in the set
$\{q_{1,1},q_{2,1}\}$ by the symbol $q_{*,1}$.\vspace{0.1cm}

Every non-singular point $o\in V$ is assumed to satisfy the
regularity condition\vspace{0.1cm}

(R1) the polynomials $q_{i,j}$ in the standard order with the last
two of them removed form a regular sequence in ${\cal
O}_{o,{\mathbb P}}$.\vspace{0.1cm}

Every quadratic point $o\in V$ is assumed to satisfy a number of
regularity conditions. Note that in this case at least one of the
linear forms $q_{1,1},q_{2,1}$ is non-zero and the other one is
proportional to it. We denote a non-zero form in the set
$\{q_{1,1},q_{2,1}\}$ by the symbol $q_{*,1}$.\vspace{0.1cm}

(R2.1) The rank of the quadratic point $o\in V$ is at least
9.\vspace{0.1cm}

{\bf Remark 0.2.} When we cut $V$ by a general linear subspace
$P\subset{\mathbb P}$ of dimension 10, containing the point $o$,
we get a complete intersection $V_P\subset P\cong{\mathbb P}^{10}$
of dimension 8 with the point $o$ an isolated singularity resolved
by one blow up  $V^+_P\to V_P$, the exceptional divisor of which,
$Q_P$, is a non-singular 7-dimensional quadric.\vspace{0.1cm}

Apart from (R2.1), the quadratic point $o$ is assumed to satisfy
the condition

(R2.2) the polynomials
$$
q_{*,1},q_{1,2},q_{2,2},\dots,q_{2,d_2}
$$
in the standard order with $q_{2,d_2}$ removed, form a regular
sequence in ${\cal O}_{o,{\mathbb P}}$.\vspace{0.1cm}

Now let us consider the biquadratic points, that is, the points
$o\in V$ for which $q_{1,1}\equiv q_{2,1}\equiv 0$.\vspace{0.1cm}

(R3.1) For a general linear subspace $P\subset{\mathbb P}$ of
dimension 12, containing the point $o$, the intersection
$V_P=V\cap P$ is a complete intersection of codimension 2 in
$P={\mathbb P}^{12}$ with the point $o\in V_P$ an isolated
singularity resolved by one blow up $V^+_P\to V_P$ with the
exceptional divisor $Q_P$ which is a non-singular complete
intersection of two quadrics in ${\mathbb P}^{11}$, $\mathop{\rm
dim}Q_P=9$.\vspace{0.1cm}

Apart from (R3.1), the biquadratic point $o$ is assumed to satisfy
the condition (R3.2) the polynomials
$$
q_{1,2},q_{2,2},\dots,q_{2,d_2}
$$
form a regular sequence in ${\cal O}_{o,{\mathbb
P}}$.\vspace{0.1cm}

The subset ${\cal P}_{\rm reg}$ consists of the pairs $(f_1,f_2)$
such that the conditions (R0.1- R0.3) are satisfied and the
conditions (R1), (R2.1) and (R2.2), (R3.1) and (R3.2) are
satisfied for every non-singular, quadratic and biquadratic point,
respectively.\vspace{0.3cm}


{\bf 0.3. The structure of the proof of Theorem 0.1.} By the well
known Grothendieck's theorem \cite{CL} for every pair
$(f_1,f_2)\in{\cal P}_{\rm reg}$ the variety $V(f_1,f_2)$
satisfies the conditions of part (i) of Theorem 0.1. Therefore,
Theorem 0.1 is implied by the following two claims.\vspace{0.1cm}

{\bf Theorem 0.2.} {\it The estimate
$$
\mathop{\rm codim}(({\cal P}\backslash{\cal P}_{\rm
reg})\subset{\cal P})\geq\frac12(M-9)(M-10)-1
$$
holds.}\vspace{0.1cm}

{\bf Theorem 0.3.} {\it For every pair $(f_1,f_2)\in{\cal P}_{\rm
reg}$ the variety $V=V(f_1,f_2)$ is birationally
superrigid.}\vspace{0.1cm}

The two claims are independent of each other and for that reason
will be shown separately: Theorem 0.2 in Section 3 and Theorem 0.3 in
Sections 1 and 2.\vspace{0.1cm}

In order to prove Theorem 0.3, we fix a mobile linear system
$\Sigma\subset|nH|$ on $V$, where $H$ is the class of a hyperplane
section. All we need to show is that $\Sigma$ has no maximal
singularities. (For all definitions and standard facts and
constructions of the method of maximal singularities we refer the
reader to \cite[Chapters 2 and 3]{Pukh13z}.)\vspace{0.1cm}

Therefore, we consider the following four options:

\begin{itemize}

\item $\Sigma$ has a maximal subvariety,

\item $\Sigma$ has an infinitely near maximal singularity, the
centre of which on $V$ is not contained in the singular locus
$\mathop{\rm Sing}V$,

\item $\Sigma$ has an infinitely near maximal singularity, the
centre of which on $V$ is contained in $\mathop{\rm Sing}V$ but
not in the locus of biquadratic points,

\item $\Sigma$ has an infinitely near maximal singularity, the
centre of which on $V$ is contained in the locus of biquadratic
points.

\end{itemize}

\noindent The first two options are excluded in Section 1 (this is
fairly straightforward), where we also prove a useful technical
claim strengthening the $4n^2$-inequality in the non-singular
case. The two remaining options are excluded in Section 2 (which
is much harder and requires some additional work).\vspace{0.1cm}

Theorem 0.2 is shown in Section 3, which completes the proof of
Theorem 0.1.\vspace{0.3cm}


{\bf 0.4. Historical remarks.} The first complete intersection
(which was not a hypersurface in the projective space) that was
shown to be birationally rigid was the complete intersection of a
quadric and cubic $V_{2\cdot 3}\subset{\mathbb P}^5$, see
\cite{IskPukh96,Pukh89c} and for a modern exposition \cite[Chapter
2]{Pukh13z}. Higher-dimensional complete intersections were
studied in \cite{Pukh01,Pukh14a,Pukh13a}; as a result of that
work, birational superrigidity is now proven for all non-singular
generic complete intersections of index 1 in the projective space,
except for three infinite series $2\cdot\dots\cdot 2,
2\cdot\dots\cdot 2\cdot 3$ and $2\cdot\dots\cdot 2\cdot 4$ and
finitely many particular families.\vspace{0.1cm}

Three-dimensional complete intersections of type $2\cdot 3$ with a
double point were studied in \cite{ChG}. Birational superrigidity
of one particular family (complete intersections of type $2\cdot
4$) of four-folds was proved in \cite{Ch03}. Recently a
considerable progress was made in the study of birational geometry
of weighted complete intersections and more complicated
subvarieties \cite{Okada1,Okada2,Okada3,AhmOkada}. Note that Fano
double hypersurfaces and cyclic covers
\cite{Pukh00c,Pukh09b,Ch06a} are also complete intersections of
index two in the weighted projective space. Finally, there is a
recent paper \cite{suzuki} claiming birational superrigidity of
certain families of complete intersections of index one, but it is
based on the ideas of \cite{TdF13}, which later turned out to be
faulty \cite{TdF15} and even in the corrected version some parts
are hard to follow. The classical techniques of the method of
maximal singularities remains the only reliable approach to
showing birational rigidity.


\section{Exclusion of maximal singularities. I.\\
Maximal subvarieties and non-singular points}

In this section we exclude maximal subvarieties of the mobile
linear system $\Sigma$ (Subsection 1.1) and infinitely near
maximal singularities of $\Sigma$, the centre of which is not
contained in the singular locus of $V$ (Subsection 1.2). After
that we show an improvement of the $4n^2$-inequality (Subsection
1.3), which will be used in Section 2 in the cases where the usual
$4n^2$-inequality is insufficient.\vspace{0.3cm}

{\bf 1.1. Exclusion of maximal subvarieties.} We start with the
following claim.\vspace{0.1cm}

{\bf Proposition 1.1.} {\it The linear system $\Sigma$ has no
maximal subvarieties.}\vspace{0.1cm}

{\bf Proof.} Assume that $B\subset V$ is a maximal subvariety for
$\Sigma$. Let us consider first the case $\mathop{\rm
codim}(B\subset V)=2$. For a general linear subspace
$P\subset{\mathbb P}$ of dimension 7 the intersection $V_P=V\cap
P$ is a non-singular complete intersection of codimension 2 in
${\mathbb P}^7$, hence for the numerical Chow group of classes of
cycles of codimension 2 on $V_P$ we have
$$
A^2V_P={\mathbb Z}H^2_P,
$$
where $H_P$ is the class of a hyperplane section of $V_P$. Now the
standard arguments \cite[Chapter 2, Section 2]{Pukh13z} give the
inequality
$$
\mathop{\rm mult}\nolimits_{B\cap P}\Sigma_P\leq n,
$$
where $\Sigma_P$ is the restriction of $\Sigma$ onto $V_P$, a
mobile subsystem of $|nH_P|$. Therefore, $\mathop{\rm
mult}_B\Sigma\leq n$ and $B$ is not a maximal subvariety --- a
contradiction.\vspace{0.1cm}

Now let us consider the case $\mathop{\rm codim}(B\subset V)\geq
3$, $B\not\subset\mathop{\rm Sing}V$. In this case we have the
inequality
$$
\mathop{\rm mult}\nolimits_BZ>4n^2,
$$
where $Z=(D_1\circ D_2)$ is the self-intersection of the system
$\Sigma$, $D_i\in\Sigma$ are general divisors. As $\mathop{\rm
deg}Z=n^2\mathop{\rm deg}V=n^2d_1d_2$, we use the inequality
$$
\frac{\mathop{\rm mult}_o}{\mathop{\rm deg}}Y\leq\frac{4}{d_1d_2},
$$
which holds for any smooth point $o\in V$ and any irreducible
subvariety $Y\subset V$ of codimension 2 (see Proposition 1.3
below) to obtain a contradiction. Finally, assume that
$B\subset\mathop{\rm Sing}V$. In this case $\mathop{\rm
codim}(B\subset V)\geq 10$, so that
$$
\mathop{\rm mult}\nolimits_B\Sigma>\delta n,
$$
where $\delta\geq 7$. Therefore, we have the inequality
$$
\mathop{\rm mult}\nolimits_BZ>98 n^2,
$$
which is impossible as for any singular point $o\in V$ and
subvariety $Y$ of codimension 2 the inequality
$$
\frac{\mathop{\rm mult}_o}{\mathop{\rm deg}}Y\leq\frac{9}{d_1d_2}
$$
holds, see Propositions 2.1 and 2.2.\vspace{0.1cm}

We have excluded all options for $B$.\vspace{0.1cm}

Q.E.D. for Proposition 1.1.\vspace{0.3cm}


{\bf 1.2. Exclusion of maximal singularities, the centre of which
is not contained in the singular locus.} Our next step is the
following\vspace{0.1cm}

{\bf Proposition 1.2.} {\it The centre $B$ of maximal singularity
$E$ is contained in the singular locus} $\mathop{\rm
Sing}V$.\vspace{0.1cm}

{\bf Proof.} Assume the converse: $B\not\subset\mathop{\rm
Sing}V$. Since $B$ is not a maximal subvariety of $\Sigma$, we see
that $\mathop{\rm codim}(B\subset V)\geq 3$ and the
$4n^2$-inequality holds:
\begin{equation}\label{22.02.2016.1}
\mathop{\rm mult}\nolimits_BZ>4n^2.
\end{equation}
Now let us show the opposite inequality.\vspace{0.1cm}

{\bf Proposition 1.3.} {\it For any non-singular point $o\in V$
and any irreducible subvariety $Y$ of codimension 2 the inequality
$$
\frac{\mathop{\rm mult}_o}{\mathop{\rm deg}}Y\leq\frac{4}{d_1d_2}
$$
holds.}\vspace{0.1cm}

{\bf Proof.} We consider the general case when $d_1+2\leq d_2$;
the obvious modifications for the two remaining cases $d_2=d_1+1$
and $d_2=d_1$ are left to the reader.\vspace{0.1cm}

Our proof is identical to the proof of Theorem 2.1 on birational
superrigidity of Fano complete intersections in \cite[Chapter 3,
Section 2]{Pukh13z}, except for the only point of difference: due
to the slightly weaker regularity condition (R1) for smooth
points, the procedure of constructing intersections with
hypertangent divisors has to terminate one step sooner than in the
cited argument. In other words, we use hypertangent divisors
$$
D_1, D_2, D'_3, D''_3,\dots, D'_i,D''_i, \dots, D'_{d_1-1},
D''_{d_1-1},
$$
followed by
$$
D_{d_1},\dots,D_{d_2-3}
$$
(as usual, $D_i\in\Lambda_i$ or $D'_i,D''_i\in\Lambda_i$ are
generic divisors in the $i$-th hypertangent linear system,
$\Lambda_i\subset |iH|$, $\mathop{\rm mult}_o\Lambda_i\geq i+1$),
{\it but not} $D_{d_2-2}$ as in \cite[Chapter 3, Section
2]{Pukh13z}, since the weaker regularity condition does not allow
to make that last step.\vspace{0.1cm}

Assuming that for $Y\ni o$ the claim of Proposition 1.3 does not
hold, we apply the technique of hypertangents divisors as outlined
above, and obtain an irreducible surface $S\ni o$, satisfying the
inequality
$$
\frac{\mathop{\rm mult}_o}{\mathop{\rm deg}}S\geq
\left(\frac{\mathop{\rm mult}_o}{\mathop{\rm
deg}}Y\right)\cdot\frac21\cdot\frac32\cdot\left(\frac43\cdot
\dots\cdot\frac{d_1}{d_1-1}\right)^2
\cdot\frac{d_1+1}{d_1}\cdot\dots\cdot\frac{d_2-2}{d_2-3}=
$$
$$
=\left(\frac{\mathop{\rm mult}_o}{\mathop{\rm deg}}Y\right) \cdot
d_1\cdot\frac{d_2-2}{3}>\frac{4(d_2-2)}{3d_2}\geq 1
$$
(the last inequality in this sequence holds as $d_2\geq 8$).
Therefore, $\mathop{\rm mult}_oS>\mathop{\rm deg}S$, which is
impossible. Proposition 1.3 is shown.

Therefore, the inequality (\ref{22.02.2016.1}) is impossible.
Proof of Proposition 1.2 is complete.\vspace{0.3cm}


{\bf 1.3. An improvement of the $4n^2$-inequality.} Let us
consider the following general situation: $X$ is a smooth affine
variety, $B\subset X$ a smooth subvariety of codimension at least
3, $\Sigma_X$ a mobile linear system on $X$ such that
$$
\mathop{\rm mult}\nolimits_B\Sigma_X=\alpha n\leq 2n
$$
for some $\alpha\in(1,2]$ and positive $n\in{\mathbb Q}$, but the
pair $\left(X,\frac{1}{n}\Sigma_X\right)$ has a non-canonical
singularity with the centre $B$. In other words, for some
birational morphism $\varphi\colon\widetilde{X}\to X$ of smooth
varieties and a $\varphi$-exceptional divisor
$E\subset\widetilde{X}$, such that $\varphi(E)=B$, the
Noether-Fano inequality
$$
\mathop{\rm ord}\nolimits_E\varphi^*\Sigma_X>na(E,X)
$$
holds. By the symbol $Z_X=(D_1\circ D_2)$ we denote the
self-intersection of the mobile linear system
$\Sigma_X$.\vspace{0.1cm}

{\bf Theorem 1.1.} {\it The following inequality holds:}
$$
\mathop{\rm mult}\nolimits_BZ_X>\frac{\alpha^2}{\alpha-1}n^2
$$

{\bf Remark 1.1.} It is easy to see that the minimum of the real
function $\frac{t^2}{t-1}$ on the interval (1,2] is attained at
$t=2$, so that the theorem improves the very well known
$4n^2$-inequality \cite[Chapter 2, Theorem 2.1]{Pukh13z}. The
proof given below is based on the idea that was first used in
\cite{Pukh89c} and later in several other papers.\vspace{0.1cm}

{\bf Proof.} We follow the arguments given in \cite[Chapter 2,
Section 2]{Pukh13z}, using the notations of the proof of the
$4n^2$-inequality given there. Repeating those arguments word for
word, we

\begin{itemize}

\item resolve the singularity $E$,

\item consider the oriented graph $\Gamma$ of the resolution,

\item divide the set of vertices of $\Gamma$ into the lower part
($\mathop{\rm codim}B_{i-1}\geq 3$) and the upper part
($\mathop{\rm codim}B_{i-1}=2$),

\item employ the technique of counting multiplicities.

\item use the optimization procedure for the quadratic function
$\sum\limits^K_{i=1}p_i\nu^2_i$

\end{itemize}

\noindent and obtain the inequality
$$
\mathop{\rm
mult}\nolimits_BZ>\frac{(2\Sigma_l+\Sigma_u)^2}{\Sigma_l(\Sigma_l+\Sigma_u)}n^2,
$$
see Subsection 2.2 in \cite[Chapter 2]{Pukh13z}. Now set
$m=\frac{1}{n^2}\mathop{\rm mult}_BZ$, so that the equality just
above can be re-written as
$$
(4-m)\Sigma^2_l+(4-m)\Sigma_l\Sigma_u+\Sigma^2_u<0.
$$
As the elementary multiplicities $\nu_i=\mathop{\rm
mult}_{B_{i-1}}\Sigma^{i-1}_X$ are non-increasing, we get the
inequalities
$$
\alpha n=\nu_1\geq\nu_2\geq\dots\geq\nu_i\geq\nu_{i+1}\geq\dots,
$$
so that the Noether-Fano inequality implies the estimate
$$
\alpha(\Sigma_l+\Sigma_u)>2\Sigma_l+\Sigma_u.
$$
As $1<\alpha\leq 2$ by assumption, we conclude that
$$
\Sigma_u>\frac{2-\alpha}{\alpha-1}\Sigma_l.
$$
Now the quadratic function $\gamma(t)=t^2+(4-m)t+(4-m)$ attains
the minimum at $t=\frac12(m-4)>0$ and is negative at $t=0$.
Therefore, if $\gamma(t_0)<0$ for some
$$
t_0>\frac{2-\alpha}{\alpha-1},
$$
then
$$
\gamma\left(\frac{2-\alpha}{\alpha-1}\right)=
\left(\frac{2-\alpha}{\alpha-1}\right)^2+(4-m)
\left(\frac{2-\alpha}{\alpha-1}\right)+(4-m)<0,
$$
which easily transforms to the required inequality
$m>\alpha^2/(\alpha-1)$. Q.E.D. for Theorem 1.1.\vspace{0.1cm}

The following elementary fact will be useful in Section 2 when
maximal singularities, the centre of which is contained in the
singular locus of $V$, are excluded.\vspace{0.1cm}

{\bf Proposition 1.4.} {\it The function of real argument
$$
\beta(t)=\frac{t^3}{t-1}
$$
is decreasing for $1<t\leq\frac32$ and increasing for
$t\geq\frac32$, so that it attains its minimum on $(1,\infty)$ at
$t=\frac32$, which is equal to $\frac{27}{4}$.}\vspace{0.1cm}

\begin{center}
\includegraphics[scale=0.8]{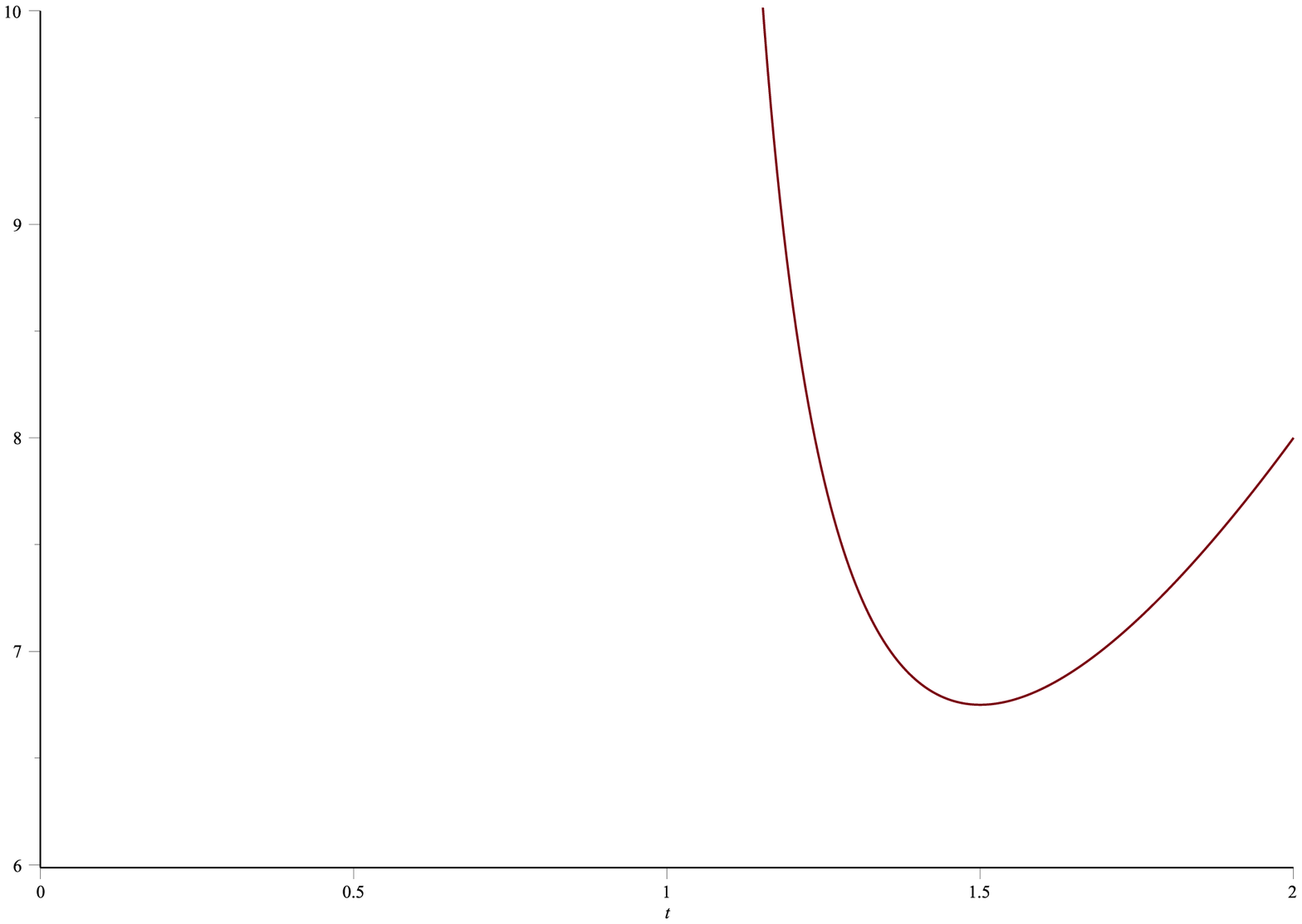}
\end{center}

{\bf Proof.} Obvious calculations. Q.E.D.


\section{Exclusion of maximal singularities. II.\\
Quadratic and biquadratic points.}

In this section we exclude infinitely near maximal singularities
of the linear system $\Sigma$, the centre of which is contained in
the singular locus of $V$. We start with using the technique of
hypertangent divisors to obtain estimates for the multiplicities
$\mathop{\rm mult}_o\Sigma$ and $\mathop{\rm mult}_oZ$, where $o$
is a general point in the centre of the maximal singularity and
$Z$ is the self-intersection of the mobile system $\Sigma$
(Subsection 2.1). After that, we consider separately the cases
when the centre is contained in the locus of the quadratic
singularities (Subsection 2.2 and 2.3) and biquadratic
singularities (Subsections 2.4 and 2.5). We make use of the
inversion of adjunction and the connectedness principle, similarly
to the arguments of Section 4 in [Book,Chapter 2], with (quite
non-trivial) modifications due to the exceptional divisor of the
blow up of the point $o$ being either a quadric or a complete
intersection of two quadrics.\vspace{0.3cm}

{\bf 2.1. The technique of hypertangent divisors.} Let
$o\in\mathop{\rm Sing}V$ be a singularity (either a quadratic or a
biquadratic point), $\sigma\colon V^+\to V$ its blow up with the
exceptional divisor $Q\subset V^+$. We consider $\sigma$ as the
resriction of the blow up $\sigma_{\mathbb P}\colon{\mathbb
P}^+\to{\mathbb P}$ of the same point $o$ on the projective space
${\mathbb P}$ with the exceptional divisor $E_{\mathbb
P}=\sigma^{-1}_{\mathbb P}(o)$, so that $Q$ is either a quadric in
a hyperplane in $E_{\mathbb P}\cong{\mathbb P}^{M+1}$ or a
complete intersection of two quadrics in $E_{\mathbb P}$. For a
generic divisor $D\in\Sigma$ set
$$
D^+\sim\sigma^*D-\nu Q
$$
for some $\nu\in{\mathbb Z}_+$; thus ${\mathop{\rm mult}_o}D=2\nu$
in the quadratic and $4\nu$ in biquadratic case. In the singular
case Proposition 1.3 has to be replaced by the following facts.
Let $Y\subset V$ be an irreducible subvariety.\vspace{0.1cm}

{\bf Proposition 2.1.} {\it Assume that $\mathop{\rm
mult}_oV$=2.\vspace{0.1cm}

{\rm (i)} If $\mathop{\rm codim} (Y\subset V)=2$, then the
inequality
$$
\frac{\mathop{\rm mult}_o}{\mathop{\rm deg}}Y\leq\frac {7}{d_1d_2}
$$
holds.\vspace{0.1cm}

{\rm (ii)} If $\mathop{\rm codim} (Y\subset V)=3$, then the
inequality
$$
\frac{\mathop{\rm mult}_o}{\mathop{\rm deg}}Y\leq\frac
{72}{7d_1d_2}
$$
holds.\vspace{0.1cm}

{\rm (iii)} The inequality $\nu\leq\sqrt{\frac72}n$
holds.}\vspace{0.1cm}

Similarly, for the biquadratic case we have\vspace{0.1cm}

{\bf Proposition 2.2.} {\it Assume that $\mathop{\rm
mult}_oV=4$.\vspace{0.1cm}

{\rm (i)} If $\mathop{\rm codim}(Y\subset V)=2$, then the
inequality
$$
\frac{\mathop{\rm mult}_o}{\mathop{\rm deg}}Y\leq\frac {9}{d_1d_2}
$$
holds.\vspace{0.1cm}

{\rm (ii)} The inequality $\nu\leq\frac32n$ holds.}\vspace{0.1cm}

{\bf Proof of Proposition 2.1.} The claim (iii) follows from (i):
for the self-intersection $Z$ of the mobile system $\Sigma$ we
have the inequality $\mathop{\rm mult}_oZ\geq 2\nu^2$. As
$\mathop{\rm deg}Z=n^2d_1d_2$, we get the inequality of part
(iii), assuming (i).\vspace{0.1cm}

In order to show the claim (i), we apply the technique of
hypertangent divisors in the same way as in the proof of
Proposition 1.3, but starting with the second hypertangent divisor
and completing the procedure with the hypertangent divisor
$D_{d_2-2}$ --- one more than in the proof of Proposition 1.3, so
that now we use the hypertangent divisors
$$
D_2,D'_3,D''_3,\dots,D'_{d_1-1},D''_{d_1-1},D_{d_1},\dots,D_{d_2-2}.
$$
If the claim (i) is not true, we obtain an irreducible surface
$S\ni o$, satisfying the inequality
$$
\frac{\mathop{\rm mult}_o}{\mathop{\rm
deg}}S\geq\left(\frac{\mathop{\rm mult}_o}{\mathop{\rm
deg}}Y\right)\cdot\frac32\cdot\left(\frac43\cdot\dots\cdot
\frac{d_1}{d_1-1}\right)^2\cdot\frac{d_1+1}{d_1}\cdot\dots\cdot
\frac{d_2-1}{d_2-2}=
$$
$$
=\left(\frac{\mathop{\rm mult}_o}{\mathop{\rm
deg}}Y\right)\cdot\frac{d_1(d_2-1)}{6}>\frac{7(d_2-1)}{6d_2}>1
$$
which is impossible. The contradiction proves the claim
(i).\vspace{0.1cm}

Finally, to show the claim (ii), we argue in exactly the same way
as above, starting with the hypertangent divisors $D'_3,D''_3$
(removing $D_2$), so that if the claim (ii) does not hold, we
obtain an irreducible surface $S\ni o$, satisfying the inequality
$$
\frac{\mathop{\rm mult}_o}{\mathop{\rm deg}}S>\frac{72}{7d_1d_2}
\cdot\left(\frac43\cdot\dots\cdot\frac{d_1}{d_1-1}\right)^2\cdot
\frac{d_1+1}{d_1}\cdot\dots\cdot\frac{d_2-1}{d_2-2}.
$$
The right hand side simplifies to
$$
\frac{72(d_2-1)}{63d_2}\geq 1
$$
for $d_2\geq 8$ which gives the desired contradiction and
completes the proof of Proposition 2.1. Q.E.D.\vspace{0.1cm}

{\bf Proof of Proposition 2.2} is very similar. First, we note
that part (i) implies part (ii) via looking at the multiplicity of
the self-intersection $Z$ at the point $o$. In order to show the
claim (i), we use the hypertangent divisors
$$
D'_3,D''_3,\dots,D'_{d_1-1},D''_{d_1-1},D_{d_1},\dots,D_{d_2-1}
$$
to obtain the required estimate. Q.E.D. for Preposition
2.2.\vspace{0.3cm}


{\bf 2.2. Exclusion of the quadratic case, part I.} In this
subsection and in the next one we assume that the centre of the
maximal singularity is contained in the singular locus
$\mathop{\rm Sing}V$ but not in the locus of biquadratic points.
We will show that this assumption leads to a contradiction. To
begin with, fix a general point $o\in V$ in the centre of the
maximal singularity.\vspace{0.1cm}

Let $\Pi\subset P$ be a general 6-plane in a 10-plane in ${\mathbb
P}$ through the point $o$. Denote by $V_{\Pi}$ and $V_P$ the
intersections $V\cap\Pi$ and $V\cap P$, respectively. By our
assumptions about the singularities of $V$, the varieties
$V_{\Pi}$ and $V_P$ are non-singular outside $o$. Let
\begin{equation}\label{02.03.2016.5}
\begin{array}{rccrcccl}
& V^+_{\Pi} & \subset & & V^+_P & \subset & V^+ &\\
\sigma_{\Pi} & \downarrow & & \sigma_P & \downarrow & &
\downarrow & \sigma\\
& V_{\Pi} & \subset & & V_P & \subset & V & \\
\end{array}
\end{equation}
be the blow ups of the point $o$ on $V_{\Pi}$, $V_P$ and $V$. The
varieties $V^+_{\Pi}$ and $V^+_P$ are non-singular. Denote the
exceptional divisors of $\sigma_{\Pi}, \sigma_P$ and $\sigma$ by
$Q_{\Pi}$, $Q_P$ and $Q$, respectively. The quadrics $Q_{\Pi}$ and
$Q_P$ are non-singular. The hyperplane sections of $V_{\Pi}$ and
$V_P$ will be written as $H_{\Pi}$ and $H_P$. Obviously, for a
general divisor $D\in\Sigma$ we have
$$
D_{\Pi}\sim nH_{\Pi}-\nu Q_{\Pi},\quad D^+_P\sim nH_P-\nu Q_P,
$$
where $D_{\Pi}=D|_{V_{\Pi}}$, $D_P=D|_{V_P}$ and the upper index +
means the strict transform. By inversion of adjunction the pairs
$(V_{\Pi},\frac{1}{n}D_{\Pi})$ and $(V_P,\frac{1}{n}D_P)$ are not
log canonical at the point $o$. As by Proposition 2.1, (iii) we
have $\nu<2n$, whereas $a(Q_{\Pi},V_{\Pi})=2$, the pair
\begin{equation}\label{01.03.2016.1}
\left(V^+_{\Pi},\frac{1}{n}D^+_{\Pi}+\frac{(\nu-2n)}{n}Q_{\Pi}\right)
\end{equation}
is not log canonical, and the centre of any of its non-log
canonical singularities is contained in the exceptional quadric
$Q_{\Pi}$ (see Lemma 4.1 in \cite[Chapter 2]{Pukh13z}). The union
of all centres of non-log canonical singularities of the pair
(\ref{01.03.2016.1}) is a connected closed set by the
Connectedness Principle \cite{Kol93,Sh93}. Therefore,

\begin{itemize}

\item either it is a point,

\item or it is a connected 1-cycle,

\item or it contains a surface on the quadric $Q_{\Pi}$.

\end{itemize}

As the union of all centres of non-log canonical singularities of
the pair (\ref{01.03.2016.1}) is a section of the union of all
centres of non-log canonical singularities of the pair
\begin{equation}\label{01.03.2016.2}
\left(V^+_P,\frac{1}{n}D^+_P+\frac{(\nu-2n)}{n}Q_P\right)
\end{equation}
by $V^+_{\Pi}\cap Q_P$ (which is a section of the non-singular
quadric $Q_P$ by a general 4-plane in $\langle Q_P\rangle$), we
see that the first option is impossible, as the smooth
7-dimensional quadric $Q_P$ can not contain a linear subspace of
dimension 4. Therefore, we conclude that the pair
(\ref{01.03.2016.2}) is not log canonical at an irreducible
subvariety $\Delta\subset Q_P$ of codimension either 1 or
2.\vspace{0.1cm}

{\bf Proposition 2.3.} {\it The case $\mathop{\rm
codim}(\Delta\subset Q_P)=1$ is impossible.}\vspace{0.1cm}

{\bf Proof.} Assume that $\Delta$ is a divisor on $Q_P$. Then by
Proposition 4.1 in \cite[Chapter 2]{Pukh13z} we have the following
estimate for the multiplicity of the self-intersection $Z_P$ of
the system $\Sigma_P=\Sigma|_{V_P}$ at the point $o$:
$$
\mathop{\rm mult}\nolimits_oZ_P\geq 2\nu^2+2\cdot
4\left(3-\frac{\nu}{n}\right)n^2
$$
(the factor 2 in the second component of the right hand side
appears since we have the inequality $\mathop{\rm deg}\Delta\geq
2$), and easy calculations give
$$
\mathop{\rm mult}\nolimits_oZ=\mathop{\rm mult}\nolimits_o Z_P\geq
16n^2,
$$
which contradicts Proposition 2.1, (i). Q.E.D. for Proposition
2.3.\vspace{0.1cm}

Therefore we assume that $\Delta\subset Q_P$ is an irreducible
subvariety of codimension 2. That option will be shown to be
impossible in the next subsection.\vspace{0.3cm}


{\bf 2.3. Exclusion of the quadratic case, part II.} Our arguments
are very similar to those in \cite[Chapter 2, Section 4]{Pukh13z}.
Let $D_1,D_2\in\Sigma$ be general divisors, $Z=(D_1\circ D_2)$ the
self-intersection of the system $\Sigma$. We can write
$$
((D_1|_{V_P})^+\circ(D_2|_{V_P})^+)=Z^+_P+Z_{P,Q}
$$
where $Z_{P,Q}$ is an effective divisor on the quadric $Q_P$. By
the standard rules of the intersection theory,
$$
\mathop{\rm mult}\nolimits_oZ=\mathop{\rm mult}\nolimits_oZ_P=
\mathop{\rm deg}(Z^+_P\circ Q_P)=2\nu^2+\mathop{\rm deg}Z_{P,Q}.
$$
Let us consider the cases $\mathop{\rm deg}\Delta=2$ (when
$\Delta$ is a section of $Q_P$ by a linear subspace of codimension
2 in $\langle Q_P\rangle$) and $\mathop{\rm deg}\Delta\geq 4$
separately. Set $\alpha=\frac{\nu}{n}<2$. Note that since
$\mathop{\rm mult}_{\Delta}\Sigma^+_P>n$ and
$\Sigma^+_P|_{Q_P}\sim\nu H_Q$, where $H_Q$ is the hyperplane
section of the quadric $Q_P$, we have the inequality $\nu>n$, so
that $\alpha>1$. By Theorem 1.1,
$$
\mathop{\rm mult}\nolimits_{\Delta}(Z^+_P+Z_{P,Q})>
\frac{\alpha^2}{\alpha-1}n^2.
$$
Assume now that $\mathop{\rm deg}\Delta\geq 4$. By Proposition
2.1, (i) we have:
$$
4\mathop{\rm mult}\nolimits_{\Delta}Z^+_P \leq\mathop{\rm
deg}(Z^+_P\circ Q_P)\leq 7n^2,
$$
so that
$$
\mathop{\rm mult}\nolimits_{\Delta}Z_{P,Q}
>\left(\frac{\alpha^2}{\alpha-1}-\frac74\right) n^2.
$$
However, for $l\in{\mathbb Z}_+$ defined by the equivalence
$$
Z_{P,Q}\sim lH_Q
$$
we have the estimate $l\geq\mathop{\rm mult}_{\Delta}Z_{P,Q}$, so
that
$$
\mathop{\rm mult}\nolimits_oZ=2(\nu^2+l)>2\left(\alpha^2+
\frac{\alpha^2}{\alpha-1}\right)-\frac74 n^2.
$$
The right hand side simplifies as
$$
2\left(\frac{\alpha^3}{\alpha-1}-\frac74\right)n^2\geq 10n^2
$$
by Proposition 1.4. Therefore, we obtained the inequality
$\mathop{\rm mult}_oZ>10n^2$, which contradicts Proposition 2.1,
(i). The case $\mathop{\rm deg}\Delta\geq 4$ is now
excluded.\vspace{0.1cm}

From now on, and until the end of this subsection, we assume that
$\mathop{\rm deg}\Delta=2$, that is, $\Delta$ is cut out on $Q_P$
by a linear subspace in $\langle Q_P\rangle$ of codimension 2. By
construction, that means that there is a subvariety
$\Delta_V\subset Q$ of codimension 2 and degree 2 (that is,
$\Delta_V$ is cut out on the quadric $Q$ by a linear subspace in
$\langle Q\rangle$ of codimension 2), such that pair
$$
\left(V^+,\frac{1}{n}\Sigma^++\frac{\nu-2n}{n}Q\right)
$$
is not log canonical at $\Delta_V$ and
$$
\Delta=\Delta_V\cap V^+_P.
$$
Let $R$ be a general hyperplane section of $V$, such that $R\ni o$
and the strict transform $R^+$ contains $\Delta_V$. Let
$Z_R=(Z\circ R)$ be the self-intersection of the mobile system
$\Sigma_R=\Sigma|_R$. Obviously,
$$
\mathop{\rm mult}\nolimits_oZ_R=\mathop{\rm mult}\nolimits_oZ +
2\mathop{\rm mult}\nolimits_{\Delta_V}Z^+.
$$
Now set $Z_{P,R}=(Z_P\circ Z_R)$. By generality of both $P$ and
$R$ we have the equalities
$$
\mathop{\rm mult}\nolimits_oZ_{P,R}=\mathop{\rm
mult}\nolimits_oZ_R, \quad \mathop{\rm
mult}\nolimits_{\Delta}Z^+_P= \mathop{\rm
mult}\nolimits_{\Delta_V}Z^+.
$$
Applying Proposition 2.1, (iii) and taking into account the
equalities above, we get the estimate
\begin{equation}\label{02.03.2016.1}
\mathop{\rm mult}\nolimits_oZ_P+2\mathop{\rm
mult}\nolimits_{\Delta}Z^+_P\leq\frac{72}{7}n^2.
\end{equation}
On the other hand, $Q_P$ is a non-singular (quadric) hypersurface,
so that by \cite[Chapter 2, Proposition 2.3]{Pukh13z} we have the
estimate
$$
\mathop{\rm deg}\nolimits_{P,Q}Z\geq 2\mathop{\rm
mult}\nolimits_{\Delta}Z_{P,Q}
$$
and for that reason
$$
\mathop{\rm mult}\nolimits_oZ_P\geq 2\nu^2+2\mathop{\rm
mult}\nolimits_{\Delta}Z_{P,Q},
$$
so that by (\ref{02.03.2016.1}) we get:
$$
\begin{array}{ccl}
\frac{72}{7}n^2 & \geq &  2\nu^2+2(\mathop{\rm
mult}\nolimits_{\Delta}Z_{P,Q}+\mathop{\rm
mult}\nolimits_{\Delta}Z^+_P)\\ & & \\
  & > &  2\left(\alpha^2+\frac{\alpha^2}{\alpha-1}\right)n^2=
2\frac{\alpha^3}{\alpha-1}n^2.
\end{array}
$$
Now we apply Proposition 1.4 and obtain the inequality
$\frac{72}{7}>\frac{27}{2}$, which is false. This contradiction
excludes the quadratic case completely.\vspace{0.3cm}


{\bf 2.4. Exclusion of the biquadratic case, part I.} In this
section and in the next one we assume that the centre of the
maximal singularity is contained in the locus of biquadratic
points. Again, we show that this assumption leads to a
contradiction. For a start, we fix a general point $o\in V$ in the
centre of the maximal singularity.\vspace{0.1cm}

Now we take a general 7-plane $\Pi$ through the point $o$ and a
general 12-plane $P\supset\Pi$. The notations $V_{\Pi}$, $V_P$
etc. have the same meaning as in quadratic case (Subsection 2.2),
the same applies to the diagram (\ref{02.03.2016.5}) and the
subsequent introductory arguments. The only difference is that the
exceptional divisors $Q_{\Pi}$ and $Q_P$ of the blow ups of the
point $o$ on $V_{\Pi}$ on $V_P$ are now non-singular complete
intersections of two quadrics. Instead of Proposition 2.1, we use
Proposition 2.2, (ii) to obtain the inequality $\nu\leq\frac32n<
2n$ and, once again, to conclude that the pair
(\ref{01.03.2016.1}) is non-log canonical. Repeating the arguments
of Subsection 2.2, we obtain the following four options for the
union of all centres of non-log canonical singularities of the
pair (\ref{01.03.2016.1}) in the biquadratic case:\vspace{0.1cm}

\begin{itemize}

\item either it is a point,

\item or it is a connected 1-cycle,

\item or it is a connected closed set of dimension 2,

\item or it contains a divisor on the 4-dimensional complete
intersection $Q_{\Pi}$.

\end{itemize}

Passing over to the pair (\ref{01.03.2016.2}) in exactly the same
way as we did it in the quadratic case, we see that the first
option is impossible as a non-singular 9-fold $Q_P$ can not
contain a linear subspace of dimension 5. Therefore, the pair
(\ref{01.03.2016.2}) is not log canonical at an irreducible
subvariety $\Delta\subset Q_P$ of codimension 1,2 or 3. The
divisorial case  ($\mathop{\rm codim}(\Delta\subset Q_P)=1$) is
excluded by the arguments of the proof of Proposition 2.3
--- in fact, we get a stronger estimate in this case:
$$
\mathop{\rm mult}\nolimits_oZ_P\geq 4\nu^2+4\cdot
4\left(3-\frac{\nu}{n}\right)n^2
$$
(as $\mathop{\rm mult}_oV_P=4$ and $\mathop{\rm deg}\Delta\geq
4$), so that
$$
\mathop{\rm mult}\nolimits_oZ=\mathop{\rm mult}\nolimits_oZ_P \geq
32n^2,
$$
which contradicts Proposition 2.2, (i).\vspace{0.1cm}

The case $\mathop{\rm codim}(\Delta\subset Q_P)=2$ is excluded by
the arguments of Subsection 2.3 as $\mathop{\rm deg}\Delta\geq 4$
and the resulting estimate $\mathop{\rm mult}_oZ> 10n^2$
contradicts Proposition 2.2, (i).\vspace{0.1cm}

It remains to exclude the last option, when $\mathop{\rm
codim}(\Delta\subset Q_P)=3$, for which there is no analog in the
quadratic case.\vspace{0.3cm}


{\bf 2.5. Exclusion of the biquadratic case, part II.} From now
on, and until the end of this section, $\Delta\subset Q_P$ is an
irreducible subvariety of codimension 3. Slightly abusing our
notations, which should not generate any misunderstanding, we show
first the following claim.\vspace{0.1cm}

{\bf Proposition 2.4.} {\it Let $Q=G_1\cap G_2\subset {\mathbb
P}^N$, $N\geq 11$, be a non-singular complete intersection of two
quadrics $G_1$ and $G_2$, $W\subset Q$ an irreducible subvariety
of codimension 2 and $\Delta \subset Q$ an irreducible subvariety
of codimension 3. Let $l\in{\mathbb Z}_+$ be defined by the
relation
$$
W\sim lH^2_Q,
$$
where $H_Q$ is the class of a hyperplane section of $Q$. Then the
inequality
$$
\mathop{\rm mult}\nolimits_{\Delta}W\leq l
$$
holds.}\vspace{0.1cm}

{\bf Proof.} Assume the converse. For a point $p\in Q$ we denote
by the symbol $|H_Q-2p|$ the pencil of tangent hyperplane sections
at that point.\vspace{0.1cm}

{\bf Lemma 2.1.} {\it Let $Y$ be an irreducible subvariety of
codimension 2, containing the subvariety $\Delta$. For a general
point $p\in\Delta$ and any divisor $T\in|H_Q-2p|$ we have}
$Y\not\subset T$.\vspace{0.1cm}

{\bf Proof of the lemma.} Assume the converse. Then for general
points $p,q\in\Delta$ and some hyperplane sections
$T_p\in|H_Q-2p|$ and $T_q\in|H_Q-2q|$ we have $Y\subset T_p\cap
T_q$, so that $Y=T_p\cap T_q$ is a section of $Q$ by a linear
subspace of codimension 2. Since $\mathop{\rm Sing}(T_p\cap T_q)$
is at most 1-dimensional (see, for instance, [Pukh00a]) and
$\mathop{\rm codim}(\Delta\subset Q)=3$, we obtain a
contradiction, varying the points $p,q$. Q.E.D. for the
lemma.\vspace{0.1cm}

We conclude that for a general point $p\in\Delta$ and an {\it
arbitrary} hyperplane section $T_p\in|H_Q-2p|$ the cycle
$W_p=(W\circ T_p)$ is well defined. It is an effective cycle of
codimension 3 on $Q$ and 2 on $T_p$ (the latter variety is a
complete intersection of two quadrics in ${\mathbb P}^{N-1}$ with
at most 0-dimensional singularities). Let $H_p\in\mathop{\rm
Pic}T_p$ be the class of a hyperplane section. Then we can write
$W_p\sim lH^2_p$. Set
$$
\Delta_p=\Delta\cap T_p.
$$
Obviously, for a general point $p$ the closed set $\Delta_p$ is of
codimension 3 on $T_p$. For any point $q\in\Delta_p$ the
inequality
$$
\mathop{\rm mult}\nolimits_qW_p>l
$$
holds. Besides, by construction $\mathop{\rm
mult}_pW_p>2l$.\vspace{0.1cm}

Now let us consider a point $q\in\Delta_p$ of general position.
Repeating the proof of Lemma 2.1 word for word (and taking into
account that the complete intersection of two quadrics $T_p$ has
zero-dimensional singularities), we see that for {\it any} divisor
$T_q\in|H_Q-2q|$ none of the components of the effective cycle
$W_p$ is contained in $T_q$, so that
$$
W_{pq}=(W_p\circ T_q)
$$
is well defined effective cycle of codimension 2 on $T_p\cap T_q$,
of codimension 3 on $T_p$ and 4 on $Q$. Since $T_q$ is an
arbitrary hyperplane section in the pencil $|H_Q-2q|$, we can
choose it to be the one containing the point $p$. Now $W_{pq}$ is
an effective cycle of codimension 6 on ${\mathbb P}^N$ of degree
$\mathop{\rm deg}W_{pq}=4$, satisfying the inequalities
$$
\mathop{\rm mult}\nolimits_pW_{pq}>2l\quad \mbox{and}\quad
\mathop{\rm mult}\nolimits_qW_{pq}>2l.
$$
Taking a general projection onto ${\mathbb P}^{N-6}$, we conclude
that the line $[p,q]\subset{\mathbb P}^N$, joining the points $p$
and $q$, is contained in the support of the cycle $W_{pq}$.
Therefore, for {\it any} point $q\in\Delta_p$ we have
$[p,q]\subset W$ and so for {\it any} point $q\in\Delta$ we have
$[p,q]\subset W$. Since $\Delta$ is not a linear subspace in
${\mathbb P}^N$ ($Q$ cannot contain linear subspaces of dimension
$N-5$) and $\mathop{\rm dim}W=N-4$, we conclude that $\Delta$ is a
hypersurface in a linear subspace of dimension $N-4$ and $W$ is
that linear subspace, which is again impossible. The proof of
Proposition 2.4 is now complete. Q.E.D.\vspace{0.1cm}

Now coming back to the biquadratic case and using the notations of
that case, we write for general divisors $D_1,D_2\in\Sigma$:
$$
((D_1|_{V_P})^+\circ (D_2|_{V_P})^+)=Z^+_P+Z_{P,Q},
$$
where again $Z_{P,Q}$ is an effective divisor on the exceptional
divisor of the blow up $\sigma_P$ of the point $o$, which is a
non-singular complete intersection of two quadrics. Again,
\begin{equation}\label{07.03.2016.1}
\mathop{\rm mult}\nolimits_oZ=\mathop{\rm mult}\nolimits_oZ_P=
\mathop{\rm deg}(Z^+_P\circ Q_P)=4\nu^2+\mathop{\rm deg}Z_{P,Q}.
\end{equation}
We set $\alpha=\frac{\nu}{n}\leq\frac32$. By Theorem 1.1,
$$
\mathop{\rm mult}\nolimits_{\Delta}(Z^+_P\circ Q_P)+ \mathop{\rm
mult}\nolimits_{\Delta}Z_{P,Q}>\frac{\alpha^2}{\alpha-1}n^2.
$$
By Proposition 2.4,
$$
\mathop{\rm mult}\nolimits_{\Delta}(Z^+_P\circ Q_P)
\leq\frac14\mathop{\rm deg}(Z^+_P\circ Q_P)= \mathop{\rm
mult}\nolimits_oZ_P.
$$
As $\mathop{\rm deg}Q_P=4$, we also have the estimate
$$
\mathop{\rm mult}\nolimits_{\Delta}Z_{P,Q} \leq\frac14\mathop{\rm
deg}Z_{P,Q},
$$
so that
$$
\mathop{\rm mult}\nolimits_oZ_P+\mathop{\rm deg}Z_{P,Q}>
4\frac{\alpha^2}{\alpha-1}n^2.
$$
Using (\ref{07.03.2016.1}), we get finally:
$$
2\mathop{\rm mult}\nolimits_oZ>4\left(\alpha^2+
\frac{\alpha^2}{\alpha-1}\right)n^2=4\frac{\alpha^3}{\alpha-1}n^2.
$$
Applying Proposition 1.4, we conclude that
$$
\mathop{\rm mult}\nolimits_oZ>\frac{27}{2}n^2,
$$
which contradicts Proposition 2.2, (i).\vspace{0.1cm}

Proof of Theorem 0.3 is now complete.


\section{Regularity conditions}

In this section we will prove Theorem 0.2 in several steps. We
first notice that
$$
\mathop{\rm codim}(({\cal P}\backslash{\cal P}_{\rm
reg})\subset{\cal P})=\mbox{ min}_{\{\ast\in S\}}\{\mathop{\rm codim}
(({\cal P}\backslash{\cal P}_{\ast})\subset{\cal P})\},
$$
where $S= \{\mbox{(R0.1),(R0.2)}, \ldots, \mbox{(R3.2)}\}$ and
$$
\mathcal{P}_{\ast}=\{(f_1,f_2)\in \mathcal{P} \, | \mbox{ the pair
satifies the regularity condition }\ast \}.
$$
We first deal with the global conditions (R0.1-R0.3) (Subsection 3.1).
Then move onto estimating the codimension of the bad set for the
condition (R1) (that is, the set of pairs $(f_1,f_2)$ that do not satisfy
that condition) and show that the same estimates work for the conditions
(R2.2) and (R3.2) (Subsections 3.2 and 3.3). Lastly, we deal with the
conditions (R2.1) and (R3.1) to get our total estimate (Subsection 3.4).\vspace{0.3cm}


{\bf 3.1. Global conditions.} We first start by splitting the condition
(R0.1) up into two conditions. The first is the irreducibility condition
for the hypersurface $\{f_1=0\}$; the set of pairs $(f_1,f_2)$ with $f_1$
irreducible is denoted by $\mathcal{P}_{\rm irred}$. The second condition is
that the hypersurface $\{f_1=0\}$ has at most quadratic singularities of
rank at least 5; the corresponding subset of ${\cal P}$ is denoted by
$\mathcal{P}_{\rm qsing\geq 5}$.\vspace{0.1cm}

{\bf Proposition 3.1.} {\it The codimension of $\mathcal{P}\setminus
\mathcal{P}_{\rm irred}$ in  $\mathcal{P}$ is at least $\frac{M(M+3)}{2}$.}\vspace{0.1cm}

{\bf Proof.} This is independent of the choice of $f_2$, hence it reduces
to looking at $f\in \mathcal{P}_{d_1,M+3}$ such that $f=g_1\cdot g_{2}$
with $\mathop{\rm deg} g_1=a$ and $\mathop{\rm deg} g_2=d_1-a$,
$a=1,2,\ldots , d_1-1$. Then we define
$$
\mathcal{F}_i = \mathcal{P}_{i,M+3}\times \mathcal{P}_{d_1-i,M+3}.
$$
Obviously, we have
$$
\mathop{\rm dim} \mathcal{P}\setminus \mathcal{P}_{\rm irred}\leq
\mathop{\rm max} \{\mathop{\rm dim} \mathcal{F}_{i}\, |\, i=1,2,\dots, d_1-1\}.
$$
We calculate:
$$
\mbox{dim }\mathcal{F}_i=\binom{i+M+2}{M+2}+\binom{d_1-i+M+2}{M+2}.
$$
By assumption $d_1 \leq \frac{M}{2}+1$. We see that this gives the maximum
dimension occuring at $i=1$, or $i=d_1-1$ as $\mathcal{F}_i=\mathcal{F}_{d_1-i}$. Then
$$
\mbox{dim }\mathcal{F}_1=(M+3)+\binom{d_1+M+1}{M+2},
$$
which immediately estimates the codimension of $\mathcal{P}\setminus
\mathcal{P}_{\rm irred}$ in $\mathcal{P}$ from below by
$$
\binom{d_1+M+2}{M+2}-\left( (M+3)+\binom{d_1+M+1}{M+2}\right)=\binom{d_1+M+1}{M+1}-(M+3).
$$
The mininal value occurs at $d_1=2$ to get the estimate claimed by our proposition.
 Q.E.D. for Proposition 3.1. \vspace{0.1cm}

{\bf Proposition 3.2.} {\it The codimension of $\mathcal{P}\setminus
\mathcal{P}_{\rm qsing \geq 5}$ in $\mathcal{P}$ is at least $\binom{M-1}{2}+1$.}\vspace{0.1cm}

{\bf Proof.} This is essentially a calculation about the rank of quadratic
forms which has been done in many places, see \cite{EP}. Q.E.D.\vspace{0.1cm}

As $\mathcal{P}_{(R0.1)}=\mathcal{P}_{\rm irred}\cap \mathcal{P}_{\rm qsing \geq 5}$,
we get that the codimension of $\mathcal{P}\setminus \mathcal{P}_{(R0.1)}$ in
$\mathcal{P}$ is at least $\binom{M-1}{2}+1$.\vspace{0.1cm}

Now we consider $\mathcal{P}_{(R0.2)}\subset \mathcal{P}$ consisting of pairs
$(f_1,f_2)$ satisfying the regularity condition (R0.2). We have two cases to
consider: the first is if the hypersurfaces contains a common component; the
second is if the intersection is non-reduced or reducible. The second case is
the only one which needs considering as the first one gives a much higher
codimension of the bad set. Fixing $f_1$ we consider the set $\mathcal{H}\subset
\mathcal{P}_{d_2,M+3}$ such that $F_1 \cap F_2$ is reducible or non-reduced.\vspace{0.1cm}

{\bf Proposition 3.3.} {\it The codimension of $\mathcal{H}$ in $\mathcal{P}_{d_2,M+3}$
is at least $\binom{M+2}{2}-2$.}\vspace{0.1cm}

{\bf Proof.} Taking into account Remark 0.1, we see that if $f_2 \in \mathcal{H}$, then:
$$
f_2|_{F_1}\in \mathcal{P}_{i,M+3}|_{F_1}\times \mathcal{P}_{d_2-i,M+3}|_{F_1},
$$
for some $i=1,2,\ldots d_2-1$. Arguing like in the proof of Proposition 3.1,
we get: the codimension of $\mathcal{H}$ in $\mathcal{P}_{d_2,M+3}$ is greater or equal than
$$
\binom{d_2+M+2}{d_2}-\left( (M+3)+\binom{d_2+M+1}{d_2-1}+\binom{d_2-d_1+M+2}{d_2-d_1} \right)
$$
$$
=\frac{1}{(M+2)!}\left( \frac{(M+2)(d_2+M+1)!}{d_2}-\frac{(d_2-d_1+M+2)!}{(d_2-d_1)!} \right)-(M+3).
$$
Using the substitution $s=d_2-d_1$, we see that for a fixed $s$ the minimum of the
above expression occurs for $d_2=s+2$ and is equal to
$$
\frac{1}{(M+2)!}\left( \frac{(M+2)(s+M+3)!}{d_2}-\frac{(s+M+2)!}{s!} \right)-(M+3).
$$
An easy check shows that this is an increasing function of $s$, so that the minimum
occurs at $s=0$ to give us the required estimate. Q.E.D. for Proposition 3.3.\vspace{0.3cm}


{\bf 3.2. Regularity conditions for smooth points.} Recall that a smooth point
satisfies the regularity condition (R1) if the homogeneous components $q_{i,j}$ in the
standard order with the last two terms (that is, the two terms of highest degree)
removed, form a regular sequence. If $d_1<d_2$, then we need
$$
W=\{q_{1,1}=q_{1,2}=\ldots =q_{1,d_1}=q_{2,1}=q_{2,2}=\ldots =q_{2,d_2-2}=0\}
$$
to be a finite set of surfaces in $\mathbb{A}^{M+2}$. If $d_1=d_2$, then we need
$$
W=\{q_{1,1}=q_{1,2}=\ldots =q_{1,d_1-1}=q_{2,1}=q_{2,2}=\ldots =q_{2,d_2-1}=0\}
$$
to be a finite set of surfaces in $\mathbb{A}^{M+2}$.\vspace{0.1cm}

 The linear forms $q_{1,1}$ and $q_{2,1}$ define the tangent space $T_pV$ at the
 point $p$, so in the case $d_2 > d_1$
$$
W=\{q_{1,2}|_{T_pV}=\ldots =q_{1,d_1}|_{T_pV}=q_{2,2}|_{T_pV}=\ldots
=q_{2,d_2-2}|_{T_pV}=0\}\subset \mathbb{A}^M
$$
and similarly for the case $d_1=d_2$. Finally as all the terms above are
homogeneous we can consider the projective variety defined by the same
equations in the projectivized tangent space. Denote this by $\widetilde{W}
\subset \mathbb{P}^{M-1}$. We have now redefined the regularity condition
under consideration to be $\mathop{\rm codim} (\widetilde{W}\subset
\mathbb{P}^{M-1}) = M-2$, that is, $\widetilde{W}$ is a finite set of curves.\vspace{0.1cm}

{\bf Proposition 3.4.} {\it The codimension of $\mathcal{P}\setminus \mathcal{P}_{(R1)}$
in $\mathcal{P}$ is at least}
$$
\lambda(M)=\frac{(M-5)(M-6)}{2}-(M+1).
$$

{\bf Proof.} We follow the methods given in \cite{Pukh01,Pukh98b} to estimate the
codimension of the space of varieties which violate the regularity conditions. The
scheme of these methods will be briefly outlined here, firstly we introduce the
necessary definitions.\vspace{0.1cm}

We say a sequence of polynomials $p_1,p_2,\ldots p_l$ is {\it $k$-regular}, with
$k \leq l$ if the subsequence $p_1,p_2, \ldots p_k$ is regular.\vspace{0.1cm}

We re-label our polynomials in their standard ordering by $h_1=q_{1,2}$, $h_2=q_{2,2}$, etc.
Also define deg $h_i=m_i$ to get our sequence $h_1, \ldots h_{M-2}$, with
$m_i\leq m_{i+1}$ in the space
$$
\mathcal{L}=\prod_{i=1}^{M-2}\mathcal{P}_{m_i,M}.
$$
We further look at the partial products defined by:
$$
\mathcal{L}_k=\prod_{i=1}^{k}\mathcal{P}_{m_i,M}.
$$
We also define
$$
Y_k(p)=\{(h_{\ast})\in \mathcal{L}_k \mbox{ } | \mbox{ } (h_{\ast}) \mbox{ is
a nonregular sequence at the point } p\},
$$
emphasising the choice of fixing the point $p$ as our origin of affine coordinates.
We will now consider $k=1,2,\ldots, M-2$ and denote
$$
Y(p)=\bigcup_{k=1}^{M-2}Y_k(p),
$$
the set of sequences which are not regular at some stage. Clearly, it is sufficient
to check that the codimension of $Y_k$ in ${\cal L}_k$ is at least $\lambda(M)+M$.
Now we outline the two methods of estimating the codimension of the bad set, with
the most important cases considered explicitly.\vspace{0.1cm}

{\bf Method 1.} We will use this method to get estimates for all cases but the one
 when the regularity fails at the last stage, this method is given in \cite{Pukh98b}.\vspace{0.1cm}

{\bf Case 1.}  For a start, let us consider the trivial case $k=1$. Here
$$
Y_1(p)=\{h_1\equiv 0\in \mathcal{P}_{2,M}\},
$$
so that
$$
\mathop{\rm codim} (Y_1(x)\subset {\mathcal{L}_1})=
\mathop{\rm dim}\mathcal{P}_{2,M}=\binom{M+1}{2}.
$$

{\bf Case 2.} Now assume that $k=2$. This is the first non-trivial case
and all the following cases follow this method. We have that
$$
Y_2(p)=\{(h_1,h_2)\in \mathcal{P}_{2,M}\times \mathcal{P}_{2,M}\mbox{ }
|\mbox{ } \mbox{codim}\{h_1=h_2=0\}<2\}.
$$
Now we have $Q=\{h_1=0\}=\bigcup Q_i \subset \mathbb{P}^{M-1}$, the
decomposition into its irreducible components and we assume that $h_1\not\equiv 0$.
Pick a general point $r\in \mathbb{P}^{M-1}$ not on $Q_i$ and consider
the projection from this point to get the map $\pi\colon
\mathbb{P}^{M-1}\dashrightarrow \mathbb{P}^{M-2}$, so that restricting
this projection onto each $Q_i$ we get a finite map $\pi_{Q_i}$, see the
figure below.
\begin{center}
\includegraphics[scale=0.6]{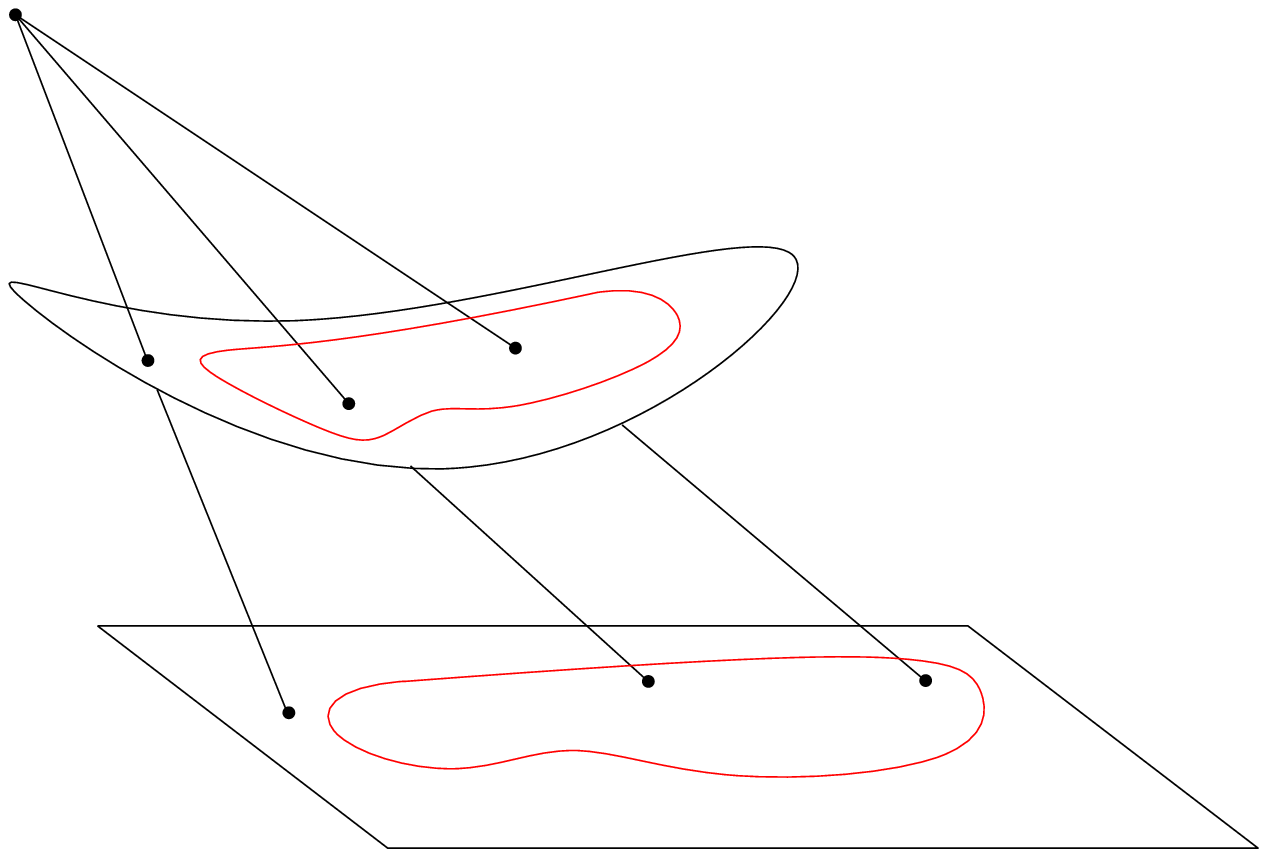}
\end{center}
Now take some $g\in H^0(\mathbb{P}^{M-2},\mathcal{O}_{\mathbb{P}^{M-2}}(2))$
and look at $\pi_{Q_i}^{\ast}(g)$: as the map is finite, we get that
$\pi_{Q_i}^{\ast}$ is injective. Therefore, for the closed subset
$$
W_2=\pi^*H^0(\mathbb{P}^{M-2},\mathcal{O}_{\mathbb{P}^{M-2}}(2))\subset
\mathcal{P}_{2,M-1}
$$
we have $W_2\cap Y_2(x)=\{0\}$. Now we know $\mathop{\rm dim} W_2=\binom{M}{2}$
so that $\mathop{\rm codim} Y_2(x)\geq \binom{M}{2}$. Therefore in the case
$k=2$ we obtain the estimate
$$
\mathop{\rm codim} (Y_2(p)\subset {\mathcal{L}_2})\geq \binom{M}{2}.
$$

{\bf The remaining cases.} We follow this method for the other values of
$k=3,\dots,M-3$; we deal with the case $k=M-2$ separately (and by means of a
different technique) later. Using this method we obtain for $k \geq 2$
($k=1$ is a special case) the inequality
$$
\mathop{\rm codim} (Y_k(p)\subset {\mathcal{L}_k})\geq \binom{\alpha_k}{\beta_k},
$$
where the values of $\alpha_k$ and $\beta_k$ are listed in the following
table ($k$ is changing from $1$ to $k=M-3$:
\[
\begin{matrix}
\alpha_k:  & M+1,& M,& M,& M-1, & M-1, & \cdots & d_2, & d_2, & d_2 , & \cdots & d_2; \\
\beta_k:  & 2, & 2, & 3, & 3, & 4, & \cdots & d_1, & d_1+1, & d_1+2, & \cdots & d_2-3.
\end{matrix}
\]
If $d_1=2$, then the smallest estimate is given by $\binom{M}{2}$, so we assume
$d_1 \geq 3$ and the smallest estimate is given by $\binom{d_2}{3}$. Now
as $d_2 \geq \frac{M}{2}+1$ we get
$$
\binom{d_2}{3} \geq \frac{M(M+2)(M-2)}{48},
$$
which is better than what we need.\vspace{0.1cm}

{\bf Method 2.} It remains to consider the case $k=M-2$. The previous projection
method outlined above in this case does not produce the estimate we need
and so we use a different method that was developed in \cite{Pukh01}. We
fix $Y^{\ast}=Y_{M-2}(p)$. Note that for any $(h_*)\in Y^*$ the sequence
$h_1,\dots,h_{M-3}$ is regular.\vspace{0.1cm}

If a sequence $(h_{\ast})$ belongs to $Y^{\ast}$ this means there exists
an irreducible component $B \subseteq Z(h_1, \ldots ,h_{M-3})$ which is
a surface with $h_{M-2}|_{B}\equiv 0$, where $Z(h_1,\ldots ,h_{M-3})\subset
\mathbb{P}^{M-1}$  is the set of common zeros of these polynomials
restricted to the projectivized tangent space.\vspace{0.1cm}

We look at the linear span $\langle B\rangle$ of $B$ and consider all
possible values of:
$$
b=\mathop{\rm codim} (\langle B\rangle\subset {\mathbb{P}^{M-1}}).
$$
Now we split $Y^{\ast}$ up into the union
$$
Y^{\ast}=\bigcup_{b=0}^{M-2}Y^{\ast}(b),
$$
where $Y^{\ast}(b)$ is the set of $(M-3)$-uples $(h_*)\in Y^*$ such that
for some irreducible curve $B\subseteq Z(h_1, \ldots , h_{M-3})$ such that
$\mathop{\rm codim} \langle B\rangle=b$, the polynomial $h_{M-2}$ vanishes
on $B$.\vspace{0.1cm}

To begin with, let us consider the case $b=0$. This means that $\langle B\rangle
= \mathbb{P}^{M-1}$. Notice that non-zero linear forms in $z_1, \ldots z_M$,
the coordinates on $\mathbb{P}^{M-1}$, do not vanish on $B$. As $h_{M-2}$ has
degree $d_2-2$ or $d_2-1$, we consider the worst case with the smaller degree,
that is, the space:
$$
W=\left\{\prod_{i=1}^{d_2-2}(a_{i,1}z_1+\ldots +a_{1,M}z_M)\right\} \subset \mathcal{P}_{d_2-1,M-1}.
$$
$W$ is a closed set with $\mathop{\rm dim} W=(M-1)(d_2-2)+1$; as
$d_2 \geq \frac{M}{2}+1$ we have $\mathop{\rm dim} W\geq \frac{(M-1)(M-2)}{2}+1$.
As $Y^{\ast}(0)\cap W=\{0\}$, we have
$$
\mathop{\rm codim} Y^{\ast}(0)\geq \frac{(M-2)(M-1)}{2}+1.
$$

Now let us deal with the case $1\leq b<M-3$. We use the technique of good
sequences and associated subvarieties, developed and described in detail in
\cite{Pukh01}.\vspace{0.1cm}

Let us fix some linear subspace $P\subset \mathbb{P}^{M-1}$ of codimension $b$.
Let $Y^{\ast}(P)$ be the set of all $(M-2)$-uples $(h_{\ast})\in Y^{\ast}(b)$
such that the closed subset $Z(h_1, \ldots, h_{M-3})$ contains an irreducible
component $B$ such that $\langle B\rangle=P$ and $h_{M-2}|_B \equiv 0$. \vspace{0.1cm}

We know \cite{Pukh01} that good sequences form an open set in the space of tuples of
polynomials and that the number of associated subvarieties is bounded from above by a
constant, depending on their degrees. Therefore, we may assume that some $(M-3-b)$
polynomials from the set $(h_1|_P,\ldots,h_{M-3}|_P)$ form a good sequence and
$B$ is one of its associated subvarieties. The worst estimate corresponds to the case
when the polynomials
$$
h_{b+1}|_P, \ldots, h_{M-3}|_P
$$
of the {\it highest} possible degrees form a good sequence and $B$ is one of its
associated subvarieties, and we will assume that this is the case.\vspace{0.1cm}

So we fix the polynomials $h_{b+1}, \ldots, h_{M-3}$ and estimate the number of
independent conditions imposed on the polynomials $h_{1}, \ldots, h_b, h_{M-2}$
by the requirement that they vanish on $B$, arguing as in the case $b=0$. Subtracting
the dimension of the Grassmannian of linear subspaces of codimension $b$ in
${\mathbb P}^{M-1}$, we get the estimate
$$
\mathop{\rm codim}(Y^{\ast}(b)\subset {\mathcal{L}}) \geq
(M-1-b)\cdot\left( \sum_{j=1}^{b}\mathop{\rm deg} h_j +
\mathop{\rm deg} h_{M-2}-b\right) + 1.
$$
Denote the right hand side of this inequality by $\theta_b$.\vspace{0.1cm}

{\bf Proposition 3.5.} {\it The following inequality
\begin{equation}\label{31.03.2016.1}
\theta_b\geq \frac{(M-2)(M-1)}{2}+1
\end{equation}
holds for all $b=1,2,\ldots, M-4$.}\vspace{0.1cm}

{\bf Proof.} It is easy to check that
$$
\gamma_b=\theta_{b+1}-\theta_b=(M-2-b)(\mathop{\rm deg} h_{b+1}-1)-
\left(\sum_{j=1}^{b}\mathop{\rm deg} h_j-b+\mathop{\rm deg} h_{M-2}\right),
$$
and since for $b\geq 2(d_1-1)$ we have $\mathop{\rm deg} h_{b+1}=
\mathop{\rm deg} h_{b}+1$, for these values of $b$ the equality
$$
\gamma_b=\gamma_{b-1}+(M-2-b)-2(\mathop{\rm deg} h_{b}-1)
$$
holds. From this equality we can see that the sequence $\theta_b$, where
$b=2(d_1-1),2d_1-1,\ldots, M-4$, has one of the following three types of behaviour:

\begin{itemize}

\item either it is non-decreasing,

\item or it is first increasing for $b=2d_1-2,\dots,a$, and then decreasing,

\item or it is decreasing.

\end{itemize}

\noindent Below it is checked that $\theta_{M-4}$ satisfies the inequality
(\ref{31.03.2016.1}). Therefore, in order to show (\ref{31.03.2016.1}) for
$b=2(d_1-1),\dots, M-4$, we only need to show this inequality for $b=2(d_1-1)$,
which is a part of the computation that we start now.\vspace{0.1cm}

Assume that $b=2l$, where $l=1,\dots, d_1-1$. Here $\theta_b=\omega_1(l)$, where
$$
\omega_1(t)=(M-1-2t)(t^2+t+d_2-2)+1.
$$
It is easy to check that $\omega'_1(t)\geq 0$ for $1\leq t\leq t_1$ for some $t_1>1$, and
$\omega'_1(t)< 0$ for $t>t_1$, so that the function of real argument $\omega_1(t)$ is
first increasing (on the interval $[1,t_1]$) and then decreasing (on $[t_1,\infty)$).
It follows that
$$
\mathop{\rm min} \{\theta_{2l}\,|\, l=1,\dots, d_1-1\}=
\mathop{\rm min} \{\theta_{2},\, \theta_{2(d_1-1)}\}.
$$
Now $\theta_2=\omega_1(1)=(M-3)d_2+1\geq \frac12 (M+2)(M-3)+1$, which satisfies
(\ref{31.03.2016.1}).\vspace{0.1cm}

Let us consider the second option: for $t=d_1-1$ we get
$$
\omega_1(d_1-1)=(M-2d_1+1)(d^2_1-2d_1+M)+1.
$$
As $2d_1-2\leq M-4$, we get the bound $d_1\leq \frac{M}{2}-1$. Looking
at the derivative of the function
$$
\omega_2(t)=(M-2t+1)(t^2-2t+M)+1,
$$
we conclude that its minimum on the interval $[2,\frac{M}{2}-1]$ is attained
at one of the endpoints, so is equal to the minimum of the two numbers:
$$
M(M-3)+1\quad \mbox{and}\quad \frac34 (M^2-4M+12)+1.
$$
Clearly, both satisfy the inequality (\ref{31.03.2016.1}).\vspace{0.1cm}

In order to complete the proof of our proposition, it remains to consider
the case $B=2l+1$, where $l=0,\dots, d_1-2$. Here $\theta_b=\omega_3(l)$,
where
$$
\omega_3(t)=(M-2-2t)(t^2+2t+d_2-1)+1.
$$
For $d_1\geq 3$ it is easy to check that the function $\omega_3(t)$ behaves
similarly to $\omega_1(t)$, first increasing and then decreasing, so that
it is sufficient ti show that $\omega_3(0)$ and $\omega_3(d_1-2)$ satisfy
the estimate (\ref{31.03.2016.1}). Indeed,
$$
\omega_3(0)=(M-2)(d_2-1)+1
$$
satisfies (\ref{31.03.2016.1}) as $d_2\geq \frac{M}{2}+1$ and for $t=d_1-2$
we get $\omega_3(d_1-2)=\omega_4(d_1)$, where
$$
\omega_4(t)=(M-2t+2)(t^2-3t+M-1)
$$
and easy computations show that (\ref{31.03.2016.1}) is satisfied here as
well.\vspace{0.1cm}

Finally, in the case $d_1=2$ we get the number
$$
\omega_3(0)=(M-2)(M-1)+1.
$$

Now the only case to consider is $b=M-4$. Here we get
$$
\mathop{\rm codim} (Y^{\ast}(b)\subset {\cal L}) \geq \frac{3}{4}(M^2-4M+6)+1.
$$
Proof of Proposition 3.5 is complete. Q.E.D. \vspace{0.1cm}

In order to complete the proof of Proposition 3.4, we have to consider
the only remaining case $b=M-3$. Here $\langle B\rangle=\mathbb{P}^2$,
which clearly implies $B\subset \mathbb{P}^{M-1}$ itself is a plane.
We do an easy dimension count, for a polynomial $h$ to satisfy $h|_B \equiv 0$
with $\mathop{\rm deg} h=e$ we get a closed algebraic set of polynomials
of codimension $\binom{e+2}{2}$ in $\mathcal{P}_{e,M}$. Therefore
$$
\mathop{\rm codim} (Y^{\ast}(M-3)\subset {\cal L}) \geq
\sum_{i=1}^{M-2}\binom{m_i+2}{2}-3(M-3).
$$
The sum takes the minimum value when $d_1=d_2$ and then we have the estimate
$$
\mathop{\rm codim} (Y^{\ast}(M-3)\subset {\cal L}) \geq
\frac{M(M+4)(M+2)}{24}-3M+1.
$$
Combining the results of both methods and simple calculation gives the estimate
$$
\mathop{\rm codim} (Y(p)\subset {\cal L}) \geq \frac{(M-5)(M-6)}{2}+1.
$$
Now Proposition 3.4 follows from a standard dimension count argument.\vspace{0.1cm}

{\bf Remark 3.3.} This is clearly not the tightest bound possible; however,
in Proposition 3.8 we have a weaker estimate.\vspace{0.1cm}


{\bf 3.3. Regularity conditions for singular points.} Recall that a point
is a quadratic singularity if $q_{1,1}$ and $q_{2,1}$ are proportional and
at least one of the terms is non-zero. We say a point is a biquadratic
singularity is $q_{1,1}=q_{2,1}=0$. The regularity conditions (R2.2) and
(R3.2) for both of these cases are similar to the smooth case (R1). The
arguments used for smooth points (R1) follow in a similar way for the two
cases (R2.2) and (R3.2). For quadratic points we work in $\mathbb{P}^M$
and for biquadratic points we work in $\mathbb{P}^{M+1}$, instead of
$\mathbb{P}^{M-1}$ and calculations are almost identical. We obtain
larger estimates for the codimension of non-regular sequences given
below.\vspace{0.1cm}

{\bf Proposition 3.6.} {\it The codimension of
$\mathcal{P}\setminus \mathcal{P}_{\ast}$ in $\mathcal{P}$ is at least
$$
\lambda(M)=\frac{(M-5)(M-6)}{2}-(M+1).
$$
for $\ast = $(R2.2) and (R3.2).} \vspace{0.1cm}

{\bf Proof.} We will outline the proof for the quadratic case (R2.2)
and the biquadratic case is treated in the same way. Instead of
restricting to the tangent space we restrict to the Zariski tangent
space $\{q_{i,1}=0\}$ for which ever $q_{i,1}$ is no-zero and work
in $\mathbb{P}^M$. We now have one extra polynomial to get our
standard ordering to be given by $h_1,\ldots ,h_{M-1}$ and our
polynomials now belong to $\mathcal{P}_{m_i,M+1}$. For the method 1,
case 1 we get the estimate:
$$
\mathop{\rm codim} (Y_1(x)\subset {\mathcal{L}_1})=
\mathop{\rm dim}\mathcal{P}_{2,M+1}=\binom{M+2}{2}.
$$
The remaining cases follow in the same way with the table given now
\[
\begin{matrix}
\alpha_k:  & M+2,& M+1,& M+1, & M, & \cdots & d_2+1, & d_2+1, & d_2+1 , & \cdots & d_2+1; \\
\beta_k:  & 2, & 2, & 3, & 3, & \cdots & d_1, & d_1+1, & d_1+2, & \cdots & d_2-2.
\end{matrix}
\]
Note that we get an extra term as we have an extra polynomial $h_{M-2}$.
Again if $d_1=2$, then the minimum is given by $\binom{M+1}{2}$ and if
$d_1\geq 3$, then the minimum is given by $\binom{d_2+1}{3}$. Now when
using the method 2 for the last case $k=M-1$, we first get
$\mbox{codim }Y^{\ast}(0)\geq \frac12 M^2+1$, so that in the notations
 of the proof of Proposition 3.5 we have possible values $b=1,\ldots ,M-2$.
 For $b < M-2$ we consider good sequences and get that:
$$
\mathop{\rm codim}(Y^{\ast}(b)\subset {\mathcal{L}}) \geq
(M-b)\cdot\left( \sum_{j=1}^{b}\mathop{\rm deg} h_j +
\mathop{\rm deg} h_{M-1}-b\right) + 1.
$$
It follows easily that
$$
\mathop{\rm codim}(Y^{\ast}(b)\subset {\mathcal{L}}) \geq
(M-1-b)\cdot\left( \sum_{j=1}^{b}\mathop{\rm deg} h_j +
\mathop{\rm deg} h_{M-2}-b\right) + 1,
$$
for $b=1,\ldots ,M-3$. For $b=M-2$ we now get
$$
\mathop{\rm codim}(Y^{\ast}(M-2)\subset \mathcal{L}) \geq
\sum_{i=1}^{M-1}\binom{m_i+2}{2}-3(M-2),
$$
and again see the estimate in the case (R1) works here also.
Q.E.D. \vspace{0.1cm}

We are left with the remaining two cases to consider now, that is,
(R2.1) and (R3.1). \vspace{0.1cm}

{\bf Proposition 3.7.} {\it The codimension of the set of complete
intersections with quadratic singularities of rank at most 8, that
is, the set $\mathcal{P}\setminus \mathcal{P}_{(R2.1)}$ in
$\mathcal{P}$ is at least $\binom{M-5}{2}+1$.}\vspace{0.1cm}

{\bf Proof.} Without loss of generality assume $q_{1,1}\neq 0$
and $q_{2,1}=\lambda q_{1,1}$ with $\lambda \in \mathbb{C}$. The
rank of the quadratic point is then given by the rank of the
quadratic form $(q_{2,2}-\lambda q_{1,2})$. The result is due
now to well know results on the codimension of quadrics of rank
at most $k$ (here $k=8$), see, for instance, \cite{EP}, where a
similar computation has been done for Fano hypersurfaces. Q.E.D.
for Proposition 3.7.\vspace{0.1cm}

{\bf Proposition 3.8.} {\it The codimension of the set
violating the condition (R3.1), that is the set
$\mathcal{P}\setminus \mathcal{P}_{(R3.1)}$ in $\mathcal{P}$
is at least $\binom{M-9}{2}-1$.}\vspace{0.1cm}

{\bf Proof.} Here we work with the space
$$
{\cal Q}= {\cal P}_{2,M+2}\times {\cal P}_{2,M+2}
$$
of pairs of quadratic forms on ${\mathbb P}^{M+1}$ (the latter
projective space interpreted as the exceptional divisor of the
blow up of a point $o\in {\mathbb P}^{M+2}$). Let $(g_1,g_2)\in
{\cal Q}$ be a pair of forms. The codimension of the closed set
of quadratic forms of rank less than 5 is $\frac{(M-4)(M-3)}{2}$,
so removing a closed set of that codimension we may assume that
$\mathop{\rm rk} g_1\geq 5$. This means that the quadric
$G_1=\{g_1=0\}$ is factorial, $\mathop{\rm Pic} G_1=\mathop{\rm
Cl} G_1={\mathbb Z} H_{G_1}$, where $H_{G_1}$ is the class of a
hyperplane section. Now for $g_2|_{G_1}$ to be non-reduced or
reducible it has to split up into hyperplane sections which gives
dimension $2M+4$. This has codimension $\frac{(M+2)(M-1)}{2}$ in
${\cal P}_{2,M-2}$. Therefore, removing a closed set of codimension
$\frac{(M-4)(M-3)}{2}$, we obtain a set ${\cal Q}^*\subset {\cal Q}$
of pairs $(g_1,g_2)$ such that the closed set $\{g_1=g_2=0\}$ is an
irreducible and reduced complete intersection of codimension 2.\vspace{0.1cm}

Let us consider the singular set of such a complete intersection,
which we denote by $\mathop{\rm Sing} (g_1,g_2)$. Note that
$\mathop{\rm Sing} (g_1,g_2)$ is the set of the points $p\in \{g_1=g_2=0\}$
where the Jacobian matrix of $g_1$ and $g_2$ has linearly dependent rows,
that is, there exists some $[\lambda_1:\lambda_2]\in {\mathbb P}^1$
with $p\in \mathop{\rm Sing} \{\lambda_1g_1+\lambda_2 g_2\}$ (where the
symbol $\mathop{\rm Sing} (g)$ denotes the singular locus of the
hypersurface $\{g=0\}$). Therefore,
$$
\mathop{\rm Sing} (g_1,g_2)\subset
\mathop{\bigcup}\limits_{[\lambda_1:\lambda_2]\in {\mathbb P}^1}
\mathop{\rm Sing} \{\lambda_1g_1+\lambda_2 g_2\},
$$
so that if
\begin{equation}\label{01.04.2016.1}
\mathop{\rm codim} (\mathop{\rm Sing} (g_1,g_2)\subset
\{g_1=g_2=0\})\leq k,
\end{equation}
then the line joining $g_1$ and $g_2$ in ${\cal P}_{2,M+2}$
meets the closed set of quadratic forms of rank at most $(k+2)$. We
conclude that the set of pairs $(g_1,g_2)\in {\cal Q}^*$
satisfying the inequality (\ref{01.04.2016.1}), has codimension at
least
$$
\frac{(M-k+1)(M-k)}{2}-1
$$
in ${\cal Q}$. Putting $k=10$ (and comparing the result with the
codimension of the complement ${\cal Q}\setminus {\cal Q}^*$ obtained
at the previous step), we complete the proof. Q.E.D. for Proposition 3.8.
\vspace{0.1cm}

Now the last thing to do is to compare the codimensions of the bad sets
for all regularity conditions and to find the minimum.\vspace{0.1cm}

Proof of Theorem 0.2 is now complete.

\begin{flushleft}
Department of Mathematical Sciences,\\
The University of Liverpool\\
e-mail: {\it pukh@liverpool.ac.uk, D.Evans2@liverpool.ac.uk}
\end{flushleft}

\end{document}